\documentclass[12pt,reqno]{article}

\usepackage[usenames]{color}
\usepackage{amssymb}
\usepackage{graphicx}
\usepackage{amscd}

\usepackage{lscape}
\usepackage{tikz}
\usepackage{tikz-cd}
\usepackage{tkz-fct}
\usepackage{pgfplots}
\usetikzlibrary{matrix}
\usetikzlibrary{fit,shapes}
\usetikzlibrary{positioning, calc}
\usepackage{colortbl}

\usepackage{amsthm}
\newtheorem{theorem}{Theorem}

\theoremstyle{definition}

\newtheorem{example}[theorem]{Example}

\usepackage[colorlinks=true,
linkcolor=webgreen, filecolor=webbrown,
citecolor=webgreen]{hyperref}

\definecolor{webgreen}{rgb}{0,.5,0}
\definecolor{webbrown}{rgb}{.6,0,0}

\usepackage{color}
\usepackage{float}
\usepackage{graphics,amsmath,amssymb}
\usepackage{amsfonts}
\usepackage{latexsym}
\usepackage{epsf}

\newcommand{\seqnum}[1]{\href{http://oeis.org/#1}{\underline{#1}}}

\begin{document}

\begin{center}
\vskip 1cm{\LARGE\bf Elliptic Curves, Riordan arrays and Lattice Paths} \vskip 1cm \large
Paul Barry\\
School of Science\\
South East Technological University\\
Ireland\\
\href{mailto:pbarry@wit.ie}{\tt pbarry@wit.ie}
\end{center}
\vskip .2 in

\begin{abstract} In this note, we show that to each elliptic curve of the form $$y^2-axy-y=x^3-bx^2-cx,$$ we can associate a family of lattice paths whose step set is determined by the parameters of the elliptic curve. The enumeration of these lattice paths is by means of an associated Riordan array. The curves and the paths have associated Somos $4$ sequences which are essentially the same. For the curves the link to Somos $4$ sequences is a classical result, via the elliptic divisibility sequence. For the paths the link is via a Hankel transform. \end{abstract}

\section{Introduction}
The equation
$$y^2+a_1xy+a_3y=x^3+a_2x^2+a_4x+a_6$$ of an elliptic curve is cubic in $x$, but quadratic in $y$. Thus the usual formula for the solution of a quadratic,
$$y^2+ax+b=0,$$ namely
$$y=\frac{-b \pm \sqrt{b^2-4ac}}{2a}$$ comes into play when solving for $y$. We now note that
$$\frac{-b - \sqrt{b^2-4ac}}{2a}=-\frac{c}{b}C\left(\frac{ac}{b^2}\right),$$ where
$$C(x)=\frac{1-\sqrt{1-4x}}{2x}$$ is the generating function of the Catalan numbers $C_n=\frac{1}{n+1}\binom{2n}{n}$. The Catalan numbers occur in the enumeration of many combinatorial objects \cite{Stan1, Stan2, StanC}, leading one to explore links between elliptic curves and appropriate combinatorial objects. In this note, we shall explore links between elliptic curves \cite{Sil, Wash} and lattice paths. It will then be interesting to see how the dynamics of elliptic curves \cite{Hone, Swart}, as manifested by the group law on such curves, and the properties of their elliptic division polynomials, is reflected in the combinatorics of the corresponding lattice paths. A bridge between them, which we will explore, is the elliptic divisibility sequence of an elliptic curve at a particular point on it, and the Hankel transform \cite{Layman} of a sequence that defines the Riordan array that enumerates the corresponding lattice paths from $(0,0)$ to $(n,k)$ \cite{Path}. We shall see that the geometry of the curve at the given point is then manifested in the continued fraction form \cite{Wall} of the generating function of this sequence.

Known integer sequences that occur in this note will be referred to by their On-Line Encyclopedia of Integer Sequences (OEIS) designation \cite{SL1, SL2}. Riordan arrays \cite{Book1, Book2} are often characterized by their $A$- sequence, but in this work, the $A$-matrix characterization of Riordan arrays is more relevant \cite{Merlini}. Relevant Riordan array links to lattice paths will be found in that reference.

An $(r,s)$ Somos $4$ sequence is a sequence $a_n$ that satisfies a recurrence of the type
$$a_n = \frac{r a_{n-1}a_{n-3} + s a_{n-2}^2}{a_{n-4}},$$ once $a_0, a_1, a_2$ and $a_3$ have been specified \cite{Chang, Hone, Yura}.
\section{Motivating examples}
We begin by taking  specific elliptic curves and proceeding to associate to them
\begin{itemize}
\item Two families of lattice paths
\item Two Riordan arrays  in the Bell subgroup of the Riordan group \cite{Book1, Book2, SGWW}
\item Two Jacobi continued fractions that reflect the geometry of the curves at the chosen point P.
\end{itemize}
The $(n,k)$-th element of the Riordan arrays will count the number of paths from $(0,0)$ to $(n,k)$ for the specific family of paths in question. The lattice paths will be determined by their step set (with algebraic multiplicities). The Hankel transform of the enumerator of paths from $(0,0)$ to $(n,0)$ will coincide essentially with the elliptic divisibility sequence of the elliptic curve.
\begin{example}\label{ex}
We start by looking at the elliptic curve
$$E_1: y^2+xy-y=x^3+2x^2+x,$$ and the point $(0,0)$ on it. Solving for $y$, we obtain
$$y_1=\frac{1-x-\sqrt{1+2x+9x^2+4x^3}}{2}$$ and
$$y_2=\frac{1-x+\sqrt{1+2x+9x^2+4x^3}}{2}.$$
Considering $y_1$ and $y_2$ as the generating functions of sequences, their expansions begin, respectively,
$$0, -1, -2, 1, 3, -7, -4,\ldots,$$ and
$$1, 0, 2, -1, -3, 7, 4,\ldots.$$ At this stage, it is not clear how to proceed. Noting that taking Hankel transforms is one of our goals, and that Hankel transforms are invariant under multiplication by $(-1)^n$ of the original sequence, we see that from this point of view, both expansions above are essentially  equivalent after the second term. A bigger picture to take into account is that the curve $E_1$ has many birationally equivalent correlates. Thus for instance the following elliptic curves are all equivalent to $E_1$.
\begin{align*}
E_2: y^2+3xy-y&=x^3+2x,\\
E_3: y^2-3xy-y&=x^3-x,\\
E_4: y^2-xy-y&=x^3+2x^2.
\end{align*}
They can each be reduced to the elliptic curve
$$y^2=x^3-\frac{43}{16}x-\frac{209}{164}.$$
For each of the above curves, we tabulate the expansions of the solutions for $y$.
\begin{center}
\begin{tabular}{|c|c|c|}
\hline
Curve & Solution for $y$ & Expansion of $y$\\
\hline
$E_1$   & $y_1=\frac{1-x-\sqrt{1+2x+9x^2+4x^3}}{2}$ & $0, -1, -2, 1, 3, -7, -4,\ldots$ \\
      & $y_2=\frac{1-x+\sqrt{1+2x+9x^2+4x^3}}{2}$ & $1, 0, 2, -1, -3, 7, 4,\ldots$ \\
      \hline
$E_2$   & $y_1=\frac{1-3x-\sqrt{1+2x+9x^2+4x^3}}{2}$ & $0, -2, -2, 1, 3, -7, -4,\ldots$\\
       & $y_2=\frac{1-3x+\sqrt{1+2x+9x^2+4x^3}}{2} $ & $1, -1, 2, -1, -3, 7, 4,\ldots$\\
       \hline
$E_3$   & $y_1=\frac{1+3x-\sqrt{1+2x+9x^2+4x^3}}{2}$ & $0, 1, -2, 1, 3, -7, -4,\ldots$ \\
       & $y_2=\frac{1+3x+\sqrt{1+2x+9x^2+4x^3}}{2}$ & $1, 2, 2, -1, -3, 7, 4,\ldots$ \\
       \hline
$E_4$   & $y_1=\frac{1+x-\sqrt{1+2x+9x^2+4x^3}}{2}$ & $0, 0, -2, 1, 3, -7, -4,\ldots$\\
       & $y_2=\frac{1+x+\sqrt{1+2x+9x^2+4x^3}}{2}$ & $1, 1, 2, -1, -3, 7, 4,\ldots$. \\
       \hline
\end{tabular}
\end{center}
We see that what is common (up to alternating sign) to all solutions is the part of the expansion after the second term. Thus for $E_1$, we consider
$$\frac{y_1+x}{x^2}$$ which expands as $-2, 1, 3, -7, -4,\ldots$. We would still like a sequence that begins with $1$ (largely for combinatorial reasons), so we consider the generating function
$$\frac{1}{1-x-x^2\frac{y_1+x}{x^2}}=\frac{2}{1-3x+\sqrt{1+2x+9x^2+4x^3}}.$$ The expansion of this generating function will have a Hankel transform equal to that of the $\frac{y_1+x}{x^2}$, with a $1$ prepended to it. However, our primary goal now is to revert the power series $\frac{x}{1-x-x^2\frac{y_1+x}{x^2}}=\frac{2x}{1-3x+\sqrt{1+2x+9x^2+4x^3}}$. Thus we seek the solution $u=u(x)$ of the equation
$$\frac{2u}{1-3u+\sqrt{1+2u+9u^2+4u^3}}=x,$$ for which we have $u(0)=0$. We arrive at
$$u(x)=\frac{1+3x-\sqrt{1+6x+9x^2-4x^3-8x^4}}{2x^2}.$$
Letting $f(x)=u(x)$, we set $g(x)=f(x)/x$. We can express $g(x)$ as
$$g(x)=\frac{1+2x}{1+3x}C\left(\frac{x^3(1+2x)}{(1+3x)^2}\right),$$ where $C(x)$ is the generating function of the Catalan numbers. We now consider lattice paths in the positive quadrant with step set
$$\{(1,1), 2*(2,1), (2,-1), (-3)*(1,0)\}.$$
Here, $2*(2,1)$ means that we have two types of $(2,1)$ step (say, blue and red), while $(-3)*(1,0)$ again means that we have three types of horizontal step $(1,0)$, but in combining path possibilities at lattice paths, we take account of signs (``negative'' paths are subtractive). These paths then satisfy the recurrence
$$t_{n,k}=t_{n-1,k-1}+t_{n-2,k+1}+2t_{n-2,k-1}-3t_{n-1,k},$$ subject to the initial conditions $t_{n,k} =0$ for $n<0$ or $k<0$, and $t_{0,0}=1$. The matrix $(t_{n,k})$ which enumerates these signed paths
begins
$$\left(
\begin{array}{ccccccc}
 1 & 0 & 0 & 0 & 0 & 0 & 0 \\
 -3 & 1 & 0 & 0 & 0 & 0 & 0 \\
 9 & -4 & 1 & 0 & 0 & 0 & 0 \\
 -26 & 15 & -5 & 1 & 0 & 0 & 0 \\
 74 & -52 & 22 & -6 & 1 & 0 & 0 \\
 -207 & 173 & -87 & 30 & -7 & 1 & 0 \\
 569 & -556 & 324 & -132 & 39 & -8 & 1 \\
\end{array}
\right).$$
This is the Riordan array
$$\left(\frac{g(x)}{1+2x}, f(x)\right)=\left(\frac{1}{1+3x}C\left(\frac{x^3(1+2x)}{(1+3x)^2}\right),\frac{x(1+2x)}{1+3x}C\left(\frac{x^3(1+2x)}{(1+3x)^2}\right)\right).$$
If we now specify the additional initial condition of $t_{1,0}=-1$, then the new matrix $(t_{n,k})$ begins
$$\left(
\begin{array}{ccccccc}
 1 & 0 & 0 & 0 & 0 & 0 & 0 \\
 -1 & 1 & 0 & 0 & 0 & 0 & 0 \\
 3 & -2 & 1 & 0 & 0 & 0 & 0 \\
 -8 & 7 & -3 & 1 & 0 & 0 & 0 \\
 22 & -22 & 12 & -4 & 1 & 0 & 0 \\
 -59 & 69 & -43 18 & -5 & 1 & 0 \\
 155 & -210 & 150 & -72 & 25 & -6 & 1 \\
\end{array}
\right).$$
This is the  Riordan array $(g(x), f(x))=(g(x), xg(x))$ of Bell type \cite{Book2}. It enumerates lattice paths in the positive quadrant with the step set
$$\{(1,1), 2*(2,1), (2,-1), (-3)*(1,0)\}$$ other than for the initial horizontal step from the origin, which has multiplicity $-1$.
The link between the Riordan array and the lattice path may be seen by expressing the bivariate generating function $\frac{g(x)}{1-xyg(x)}$ of the array in the form
$$\frac{g(x)}{1-xyg(x)}=\frac{1+x(3-2y)-4x^2y-\sqrt{1+6x+9x^2-9x^3-8x^4}}{2xy(2x^2y+x^2/y+xy-3x-1)}.$$ The step set can be read from the denominator.
Thus we have gone from the elliptic curve to a set of (signed) lattice paths in the positive quadrant, counted by a Riordan array.

From the above, we see that the expansion of $g(x)$ (which was obtained by solving for $y$ and the inversion process) begins
$$ 1, -1, 3, -8, 22, -59, 155, -396, 978, -2310, 5122, -10260, 16752,\ldots.$$
Calculating its Hankel transform, we find the sequence that begins
$$1, 2, 1, -7, -16, -57, -113, 670, 3983, 23647, 140576,\ldots.$$
Now the elliptic divisibility sequence of $E_1$ at $P=(0,0)$ begins
$$0,1,1,2,1,-7, -16, -57, -113, 670, 3983, 23647, 140576,\ldots.$$
The group structure of the elliptic curve is reflected in the following continued fraction form of the generating function $g(x)$. Let
$nP=([nP]_1, [nP]_2)$ denote the point $$P+P+\cdots+P \quad (n\,\, \text{times}),$$ where $P=(0,0)$ on the curve $E_1$.
Then we have
$$g(x)=\cfrac{1}{1-\frac{[2P]_2-1}{[2P]_1}x+x-\cfrac{[2P]_1x^2}{1-\frac{[3P]_2-1}{[3P]_1}x+x-\cfrac{[3P]_1x^2}{1-\frac{[4P]_2-1}{[4P]_1}x+x-\cdots}}}.$$
In particular, the Hankel transform $h_n$ of the expansion of $g(x)$ is given by
$$h_n = \prod_{k=0}^n [(k+2)P]_1^{n-k}.$$
The generating function $g(x)=\frac{1+2x}{1+3x}C\left(\frac{x^3(1+2x)}{(1+3x)^2}\right)$ can be transformed by applying a binomial transform twice to obtain the generating function
$$g_2(x)=\frac{1}{1-2x}g\left(\frac{x}{1-2x}\right)=\frac{1}{1-x-2x^2}C\left(\frac{x^3}{(1-x-2x^2)^2}\right).$$
This generating function expands to give the sequence (\seqnum{A025243}) that begins
$$1, 1, 3, 6, 14, 33, 79, 194, 482, 1214, 3090,\ldots.$$ As we have used a binomial transform to obtain this sequence, it again will have the Hankel transform that begins
$$1, 2, 1, -7, -16, -57,\ldots.$$ This sequence counts lattice paths from $(0,0)$ to $(n,0)$ with step set
$$\{(1,-1),(1,0),2*(2,0),(2,1)\}.$$ The Bell-type Riordan array $(g_2(x), xg_2(x))$, which begins
$$\left(
\begin{array}{ccccccc}
 1 & 0 & 0 & 0 & 0 & 0 & 0 \\
 1 & 1 & 0 & 0 & 0 & 0 & 0 \\
 3 & 2 & 1 & 0 & 0 & 0 & 0 \\
 6 & 7 & 3 & 1 & 0 & 0 & 0 \\
 14 & 18 & 12 & 4 & 1 & 0 & 0 \\
 33 & 49 & 37 & 18 & 5 & 1 & 0 \\
 79 & 130 & 114 & 64 & 25 & 6 & 1 \\
\end{array}
\right),$$
counts lattice paths from $(0,0)$ to $(n,k)$ with this step set. As a second binomial transform of $g(x)$, the generating function $g_2(x)$ can be expressed as the continued fraction
$$g_2(x)=\cfrac{1}{1-\frac{[2P]_2-1}{[2P]_1}x-x-\cfrac{[2P]_1x^2}{1-\frac{[3P]_2-1}{[3P]_1}x-x-\cfrac{[3P]_1x^2}{1-\frac{[4P]_2-1}{[4P]_1}x-x-\cdots}}},$$ thus making explicit its links, and those of the lattice paths that the Riordan array $(g_2(x), xg_2(x))$ counts, with the elliptic curve $E_1$, and the multiples of $P=(0,0)$ on it.

We note that the multiplies $[nP]$ of the point $P=(0,0)$ on $E_1$ have rational coordinates
$$(0,0),(0,0),(-2,1),\left(-\frac{1}{4},\frac{9}{8}\right),(14,50),\left(\frac{16}{49},-\frac{169}{343}\right),\left(-\frac{399}{256},\frac{847}{4096}\right), \left(-\frac{1808}{3249},\frac{274576}{185193}\right),\ldots.$$
The fact that the continued fractions for $g(x)$ and $g_2(x)$ expand to give an integer sequence is another manifestation of the ``integrality'' of the dynamics involved \cite{Swart}.

Since a binomial transform of a sequence leaves the Hankel transform unchanged, we can see that for general $r \in \mathbb{Z}$, the generating function
$$g_r(x)=\cfrac{1}{1-\frac{[2P]_2-1}{[2P]_1}x-rx-\cfrac{[2P]_1x^2}{1-\frac{[3P]_2-1}{[3P]_1}x-rx-\cfrac{[3P]_1x^2}{1-\frac{[4P]_2-1}{[4P]_1}x-rx-\cdots}}}$$ expands to a sequence with the same Hankel transform as before. We can write this as
$$g_r(x)=\cfrac{1}{1-\frac{[2P]_2+r[2P]_1-1}{[2P]_1}x-\cfrac{[2P]_1x^2}{1-\frac{[3P]_2+r[3P]_1-1}{[3P]_1}x-\cfrac{[3P]_1x^2}{1-\frac{[4P]_2+r[4P]_1-1}{[4P]_1}x-\cdots}}}.$$ For each $r$, this continued fraction corresponds to a curve, parameterized by $r$, that is birationally equivalent to $E_1$. Interesting sequences emerge in this ``binomial orbit''. We see the following.
\begin{align*}
t_{n-1,k-1}: &\,1,1,1,\ldots\\
t_{n-2,k-1}: &\,2,1,0,-1,-2,-3,\ldots,2-n,\ldots\\
t_{n-2,k\,\,}: &\,0,2,2,0,-4,-10,-18,\ldots,n(3-n),\ldots\\
t_{n-2,k+1}: &\,1,1,1,\ldots\\
t_{n-1,k\,\,}: &\,-3,-1,1,3,5,9,\ldots,2n-3,\ldots\\
\end{align*}

  \begin{figure}
\begin{center}
\begin{tikzpicture}
\draw[thick] (0, 0) grid (4,4);
\node[circle, draw, minimum size=.2cm, fill=red] at (1,2){a};
\node[circle, draw, minimum size=.2cm, fill=red] at (1,3){b};
\node[circle, draw, minimum size=.2cm, fill=red] at (2,3){c};
\node[circle, draw, minimum size=.2cm, fill=red] at (3,3){d};
\node[circle, draw, minimum size=.2cm, fill=red] at (2,2){e};
\node[circle, draw, minimum size=.2cm, fill=green] at (2,1){};
\end{tikzpicture}
\end{center}
\caption{Generic path set diagram for Example \ref{ex}, for the recurrence
$t_{n,k}=at_{n-1,k-1}+bt_{n-2,k-1}+et_{n-1,k}+ct_{n-2,k}+dt_{n-2,k+1}$, for the path set $\{a*(1,1), b*(2,1),c*(2,0),e*(1,0),d*(2,-1)\}$}
\end{figure}
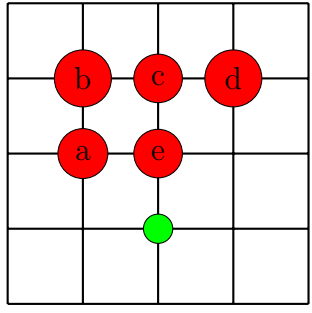

  \begin{figure}
\begin{center}
\begin{tikzpicture}
\draw[thick] (0, 0) grid (4,4);
\node[circle, draw, minimum size=.2cm, fill=red] at (1,2){1};
\node[circle, draw, minimum size=.2cm, fill=red] at (1,3){2};
\node[circle, draw, minimum size=.2cm, fill=red] at (2,3){0};
\node[circle, draw, minimum size=.2cm, fill=red] at (3,3){1};
\node[circle, draw, minimum size=.2cm, fill=red] at (2,2){-3};
\node[circle, draw, minimum size=.2cm, fill=green] at (2,1){};
\end{tikzpicture}
\end{center}
\caption{Path set for Example \ref{ex}, $g$}
\end{figure}

  \begin{figure}
\begin{center}
\begin{tikzpicture}
\draw[thick] (0, 0) grid (4,4);
\node[circle, draw, minimum size=.2cm, fill=red] at (1,2){1};
\node[circle, draw, minimum size=.2cm, fill=red] at (1,3){1};
\node[circle, draw, minimum size=.2cm, fill=red] at (2,3){2};
\node[circle, draw, minimum size=.2cm, fill=red] at (3,3){1};
\node[circle, draw, minimum size=.2cm, fill=red] at (2,2){-1};
\node[circle, draw, minimum size=.2cm, fill=green] at (2,1){};
\end{tikzpicture}
\end{center}
\caption{Path set for Example \ref{ex}, $g_1$}
\end{figure}
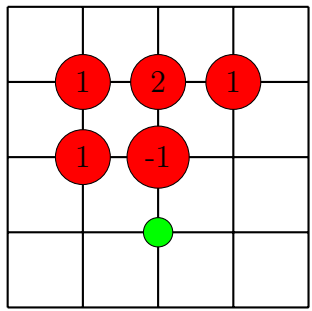

  \begin{figure}
\begin{center}
\begin{tikzpicture}
\draw[thick] (0, 0) grid (4,4);
\node[circle, draw, minimum size=.2cm, fill=red] at (1,2){1};
\node[circle, draw, minimum size=.2cm, fill=red] at (1,3){0};
\node[circle, draw, minimum size=.2cm, fill=red] at (2,3){2};
\node[circle, draw, minimum size=.2cm, fill=red] at (3,3){1};
\node[circle, draw, minimum size=.2cm, fill=red] at (2,2){1};
\node[circle, draw, minimum size=.2cm, fill=green] at (2,1){};
\end{tikzpicture}
\end{center}
\caption{Path set for Example \ref{ex}, $g_2$}
\end{figure}

  \begin{figure}
\begin{center}
\begin{tikzpicture}
\draw[thick] (0, 0) grid (4,4);
\node[circle, draw, minimum size=.2cm, fill=red] at (1,2){1};
\node[circle, draw, minimum size=.2cm, fill=red] at (1,3){-1};
\node[circle, draw, minimum size=.2cm, fill=red] at (2,3){0};
\node[circle, draw, minimum size=.2cm, fill=red] at (3,3){1};
\node[circle, draw, minimum size=.2cm, fill=red] at (2,2){3};
\node[circle, draw, minimum size=.2cm, fill=green] at (2,1){};
\end{tikzpicture}
\end{center}
\caption{Path set for Example \ref{ex}, $g_3$}
\end{figure}

  \begin{figure}
\begin{center}
\begin{tikzpicture}
\draw[thick] (0, 0) grid (4,4);
\node[circle, draw, minimum size=.2cm, fill=red] at (1,2){1};
\node[circle, draw, minimum size=.2cm, fill=red] at (1,3){-2};
\node[circle, draw, minimum size=.2cm, fill=red] at (2,3){-4};
\node[circle, draw, minimum size=.2cm, fill=red] at (3,3){1};
\node[circle, draw, minimum size=.2cm, fill=red] at (2,2){5};
\node[circle, draw, minimum size=.2cm, fill=green] at (2,1){};
\end{tikzpicture}
\end{center}
\caption{Path set for Example \ref{ex}, $g_4$}
\end{figure}

  \begin{figure}
\begin{center}
\begin{tikzpicture}
\draw[thick] (0, 0) grid (4,4);
\node[circle, draw, minimum size=.2cm, fill=red] at (1,2){1};
\node[circle, draw, minimum size=.2cm, fill=red] at (1,3){-3};
\node[circle, draw, minimum size=.2cm, fill=red] at (2,3){-10};
\node[circle, draw, minimum size=.2cm, fill=red] at (3,3){1};
\node[circle, draw, minimum size=.2cm, fill=red] at (2,2){7};
\node[circle, draw, minimum size=.2cm, fill=green] at (2,1){};
\end{tikzpicture}
\end{center}
\caption{Path set for Example \ref{ex}, $g_5$}
\end{figure}

  \begin{figure}
\begin{center}
\begin{tikzpicture}
\draw[thick] (0, 0) grid (4,4);
\node[circle, draw, minimum size=.2cm, fill=red] at (1,2){1};
\node[circle, draw, minimum size=.2cm, fill=red] at (1,3){-4};
\node[circle, draw, minimum size=.2cm, fill=red] at (2,3){-18};
\node[circle, draw, minimum size=.2cm, fill=red] at (3,3){1};
\node[circle, draw, minimum size=.2cm, fill=red] at (2,2){9};
\node[circle, draw, minimum size=.2cm, fill=green] at (2,1){};
\end{tikzpicture}
\end{center}
\caption{Path set for Example \ref{ex}, $g_6$}
\end{figure}

We end this example by looking at how the lattice path step sets change as we take successive binomial transforms of the expansion of $g(x)$. Figures 1-8 depict the steps for $g_i(x)$, where $g_i(x)$ is the $i$-th binomial transform of $g(x)$. All the resulting sequences will have the same Hankel transform, essentially the elliptic division sequence of $E_1$.
\end{example}
\begin{example}
In this next example, we look at the elliptic curve
$$y^2+2xy-y=x^3+5x^2-x,$$ and the point $(0,0)$. The elliptic divisibility sequence of this curve at $(0,0)$ begins
$$1,1,2,-9,-17,-196,593,9657,152170,\ldots.$$
We now undertake the following steps.
\begin{itemize}
\item We solve the above equation for $y$, and chose the solution $y_1=\frac{1-2x-\sqrt{1-8x+24x^2+4x^3}}{2}$ that expands as $0,1,\ldots$.
\item We form $y_2=\frac{y_1-x}{x^2}$ (discarding the first two elements of the expansion of $y_1$).
\item We form the generating function $\frac{x}{1-x-x^2y_1}$ and we revert this to get the power series $f(x)$.
\item We let $g(x)=\frac{f(x)}{x}$.
\end{itemize}
We obtain that
$$g(x)=\frac{1-2x-5x^2-\sqrt{1-4x-6x^2+16x^3+37x^4}}{2x^3}=\frac{1-3x}{1-2x-5x^2}c\left(\frac{x^3(1-3x)}{(1-2x-5x^2)^2}\right).$$
Here, $c(x)=\frac{1-\sqrt{1-4x}}{2x}$ denotes the generating function of the Catalan numbers $C_n=\frac{1}{n+1}\binom{2n}{n}$ \seqnum{A000108}.
The form of this generating function and the theory of lattice paths now tells us that the Riordan array
$(g(x), f(x))$ counts lattice from $(0,0)$ to $(n,k)$ for lattice paths whose step set is given by
$$\{(1,1), 2*(1,0),5*(2,0),(2,-1),(-3)*(2,1)\}.$$ This can be verified from the theory of lattice paths and Riordan arrays, since we have $g(x)=\frac{u}{x}$ where
$$\frac{u}{x}=1+2u+5ux+u^2x-3x.$$
The form of $g(x)$ and the theory of Hankel transforms and Somos $4$ sequences tells us that the Hankel transform of $g(x)$ is a $(1,-2)$ Somos $4$ sequence.
Calculating, we find that this Hankel transform begins
$$1, 2, -9, -17, -196, 593,\ldots.$$
The power series $g(x)$ expands to give the sequence that begins
$$1, -1, 3, 2, 17, 51, 185, 664, 2333, 8360, 29717,\ldots.$$
Taking the third binomial transform of this sequence yields the sequence that begins
$$1, -4, 18, -79, 344, -1482, 6314, -26576, 110372, -451531, 1815500,\ldots,$$ which again has a Hankel transform that begins
$$1, 2, -9, -17, -196, 593,\ldots,$$
(Hankel transforms are unchanged by binomial transforms).
The new generating function after the binomial transform can be expressed as
$$g_1(x)=\frac{1}{1+4x-2x^2}c\left(\frac{x^3}{(1+4x-2x^2)^2}\right).$$
The Riordan array $(g_1(x), x g_1(x))$ then counts lattice paths from $(0,0)$ to $(n,k)$ with step set $\{(1,1), (-4)*(1,0), 2*(2,0), (2,-1)\}$. This can be verified from the theory of lattice paths and Riordan arrays, since we have $g_1(x)=\frac{u}{x}$ where
$$\frac{u}{x}=1-4u+2ux+u^2x.$$
We note that the elements $t_{n,k}$ of the Riordan array $(g(x), f(x))$ satisfy the recurrence
$$t_{n,k}=t_{n-1,k-1}+2 t_{n-1,k}+5 t_{n-2,k}+t_{n-2,k+1}-3t_{n-2,k-1},$$ subject to the initial conditions $t_{n,k}=0$ if $n<0$ or $k<0$, $t_{0,0}=1$, and $t_{1,0}=-1$.
Similarly the elements $\tau_{n,k}$ of the Riordan array $(g_1(x), xg_1(x))$ satisfy the recurrence
$$\tau_{n,k}=\tau_{n-1,k-1}-4 \tau_{n-1,k}+2 \tau_{n-2,k}+\tau_{n-2,k+1},$$ subject to the initial conditions  $\tau_{n,k}=0$ if $n<0$ or $k<0$ and $\tau_{0,0}=1$.

We can summarize this section as follows. To the elliptic curve
$$y^2+2xy-y=x^3+5x^2-x,$$ and the point $(0,0)$ on it, we have associated
\begin{itemize}
\item The lattice paths from $(0,0)$ to $(n,k)$ with step set $\{(1,1), 2*(1,0),5*(2,0),(2,-1),(-3)*(2,1)\}$, and the lattice paths from $(0,0)$ to $(n,k)$ with step set
$\{(1,1), (-4)*(1,0), 2*(2,0), (2,-1)\}$.
\item The path enumerator Riordan arrays $(g(x), xg(x))$ and $(g_1(x), xg_1(x))$.
\end{itemize}
The Riordan array $(g(x), xg(x))$ begins
$$\left(
\begin{array}{cccccc}
 1 & 0 & 0 & 0 & 0 & 0 \\
-1 & 1 & 0 & 0 & 0 & 0 \\
 3 & -2 & 1 & 0 & 0 & 0  \\
 2 & 7 & -3 & 1 & 0 & 0  \\
 17 & -2 & 12 & -4 & 1 & 0  \\
 51 & 39 & -13 & 18 & -5 & 1  \\
\end{array}
\right).$$
The Riordan array $(g_1(x), xg_1(x))$ begins
$$\left(
\begin{array}{cccccc}
 1 & 0 & 0 & 0 & 0 & 0 \\
-4 & 1 & 0 & 0 & 0 & 0 \\
 18 & -8 & 1 & 0 & 0 & 0  \\
 -79 & 52 & -12 & 1 & 0 & 0  \\
344 & -302 & 102 & -16 & 1 & 0  \\
-1482 & 1644 & -733 & 168 & -20 & 1  \\
\end{array}
\right).$$
Note that we are considering lattice paths with algebraic multiplicities, some of which can be negative, thus entries in these matrices can be negative (paths with negative multiplicities are ``subtractive'').

In addition, the Hankel transform of both sequences $1,-1,3,2,17,\ldots$ and $1,-4,18,-79,\ldots$ is the $(1,-2)$ Somos $4$ sequence that begins
$$1, 2, -9, -17, -196, 593,\ldots.$$ This coincides essentially with the elliptic divisibility sequence of the original elliptic curve.

We note that had we started with the elliptic curve
$$y^2-2xy-y=x^3+5x^2+x,$$ then we would have arrived at the sequence
$$1, 4, 18, 81, 368, 1686, 7786, 36224, 169700,\ldots$$ with generating function
$$g_1(x)=\frac{1}{1-4x-2x^2}c\left(\frac{x^3}{(1-4x-2x^2)^2}\right).$$
The Hankel transform of this sequence begins
$$1, 2, 7, -1, -100, -351,\ldots.$$ Again, this is a $(1,-2)$ Somos $4$ sequence which essentially coincides with the elliptic divisibility sequence of the elliptic curve in question.
The Riordan array $(g_1(x), xg_1(x))$ in this case begins
$$\left(
\begin{array}{cccccc}
 1 & 0 & 0 & 0 & 0 & 0 \\
4 & 1 & 0 & 0 & 0 & 0 \\
18 & 8 & 1 & 0 & 0 & 0  \\
81 & 52 & 12 & 1 & 0 & 0  \\
368 & 306 & 102 & 16 & 1 & 0  \\
1686 & 1708 & 739 & 168 & 20 & 1  \\
\end{array}
\right),$$ and its terms satisfy the recurrence
$$t_{n,k}=t_{n-1,k-1}+4t_{n-1,k}+2 t_{n-2,k}+t_{n-2,k+1}.$$ This is the enumerator matrix for the lattice paths from $(0,0)$ to $(n,k)$ with step set $\{(1,1), 4*(1,0), 2*(2,0), (2,-1)\}$.
The generating function of this matrix can be expressed as
$$\frac{1-2x(y+2)-2x^2-\sqrt{1-8x+12x^2+12x^3+4x^4}}{2xy(xy+4x+2x^2+x^2/y-1)}.$$ We see that we can read off the elements of the step set from the denominator.
\end{example}
\section{The general case}
We now begin with the elliptic curve
$$y^2-axy-y=x^3-bx^2-cx,$$ and the point $(0,0)$ on it.
We first solve for $y$. We obtain
$$y_1=\frac{1+ax-\sqrt{1+2x(a-2c)+x^2(a^2-4b)+4x^3}}{2},$$
which expands as $0,c,-ac+b+c^2,\ldots$, or
$$y_2=\frac{1+ax+\sqrt{1+2x(a-2c)+x^2(a^2-4b)+4x^3}}{2},$$ which expands $1,a-c,ac-b-c^2,\ldots$.
Apart from signs, both expansions coincide after the second term, so we discard to first two elements to retain what is common.
Thus we consider $\frac{y_1-cx}{x^2}$, which begins $-ac+b+c^2,\ldots$.

We now take the reversion of the expression
$$\frac{x}{1-x-x^2\left(\frac{y_1-cx}{x^2}\right)}=\frac{2x}{1+x(a-2c+2)+\sqrt{1+2x(a-2c)+x^2(a^2-4b)+4x^3}},$$ to obtain the power series $f(x)=xg(x)$ where $g(x)$ is given by
\begin{scriptsize}
$$\frac{1+(a-2c+1)x}{1-x(2(c-1)-a)-x^2(a(c-1)-b-(c-1)^2)}c\left(\frac{x^3(1-x(a-2c+1))}{(1-x(2(c-1)-a)-x^2(a(c-1)-b-(c-1)^2))^2}\right).$$
\end{scriptsize}
This is the generating function of the sequence $u_n$ where
$$u_n=\sum_{k=0}^n \sum_{j=0}^{k+1}\binom{k+1}{j} c^j \sum_{i=0}^{n-3k-j}\binom{2k+i}{i}\binom{i}{n-3k-i-j}a^{2i+3k+j-n}b^{n-3k-i-j} C_k.$$
The Hankel transform of $u_n$ begins
$$1,ac-b-c^2, a^2c - a(b + 3c^2) + 2bc + 2c^3 - 1,\ldots.$$
This is a $(1, - ac + b + c^2)$ Somos $4$ sequence, which essentially coincides with the elliptic divisibility sequence of the elliptic curve.

We can show that the denominator of the  generating function of the Riordan array $(g(x), xg(x))$ contains the term
$$1-xy-x(2(c-1)-a)-x^2(a(c-1)-b-(c-1)^2)-x^2y(a-2c+1)-x^2/y,$$ thus indicating that it is linked to the family of lattice paths with step set
$$\{(1,1), (2(c-1)-a)*(1,0), (a(c-1)-b-(c-1)^2)*(2,0), (a-2c+1)*(2,1),(2,-1)\}.$$ This can be verified since the solution $u$ of the equation
$$\frac{u}{x}=1+(2(c-1)-a)u+(a(c-1)-b-(c-1)^2)ux+(a-2c+1)x+xu^2$$ is given by $u(x)=f(x)$.

We next apply the binomial transform $\left(\frac{1}{1-(a-2c+1)x}, \frac{x}{1-(a-2c+1)x}\right)$ to the generating function $g(x)$. We obtain the generating function
$$\gamma(x)=\frac{1}{1-(a-2c)x-(ac-b-c^2)x^2}c\left(\frac{x^3}{(1-(a-2c)x-(ac-b-c^2)x^2)^2}\right).$$  This is the generating function of the sequence $v_n$ given by
$$v_n=\sum_{k=0}^n \sum_{j=0}^{n-3k} \binom{2k+j}{j}\binom{j}{n-3k-j}(ac-b-c^2)^{n-3k-j}(a-2c)^{2j-n+3k}C_k.$$ As we have applied a binomial transform to arrive at this sequence, its Hankel transform will be that of $u_n$, and so remains essentially the elliptic divisibility sequence of the elliptic curve.

The denominator of the generating function $\frac{\gamma(x)}{1-xy\gamma(x)}$ of the Riordan array
$(\gamma(x), x\gamma(x)$ contains the term
$$1-xy-(a-2c)x-(ac-b-c^2)x^2+x^2/y,$$
indicating that the  Riordan array $(\gamma(x), x\gamma(x))$ is the enumerator matrix for the lattice path family with step set
$$\{(1,1), (a-2c)*(1,0), (ac-b-c^2)*(2,0), (-1)*(2,-1)\}.$$ Again, this can be verified by showing that the solution of the equation
$$\frac{u}{x}=1+(a-2c)u+(ac-b-c^2)ux-u^2x$$ is given by $u(x)=x \gamma(x)$.

  \begin{figure}
\begin{center}
\begin{tikzpicture}
\draw[thick] (0, 0) grid (4,4);
\node at (0.9,1.5){\begin{scriptsize}$1$\end{scriptsize}};
\node[circle, draw, minimum size=.2cm, fill=red] at (1,2){};
\node at (0.2,3.2){\begin{scriptsize}$a-2c+1$\end{scriptsize}};
\node[circle, draw, minimum size=.2cm, fill=red] at (1,3){};
\node  at (1.95,3.6){\begin{tiny}$a(c-1)-b-(c-1)^2$\end{tiny}};
\node[circle, draw, minimum size=.2cm, fill=red] at (2,3){};
\node at (3.2, 3.2){\begin{scriptsize}$1$\end{scriptsize}};
\node[circle, draw, minimum size=.2cm, fill=red] at (3,3){};
\node at (2.0,1.5){\begin{scriptsize}$2(c-1)-a$\end{scriptsize}};
\node[circle, draw, minimum size=.2cm, fill=red] at (2,2){};
\node[circle, draw, minimum size=.2cm, fill=green] at (2,1){};
\end{tikzpicture}
\end{center}
\caption{Path set for $g$}
\end{figure}
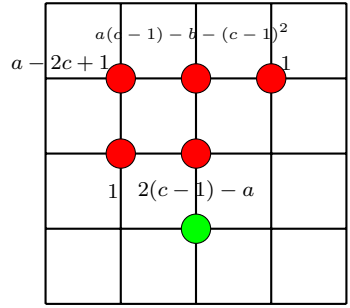

  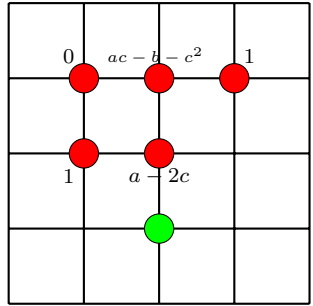
\begin{figure}
\begin{center}
\begin{tikzpicture}
\draw[thick] (0, 0) grid (4,4);
\node at (0.8,1.7){\begin{scriptsize}$1$\end{scriptsize}};
\node[circle, draw, minimum size=.2cm, fill=red] at (1,2){};
\node at (0.8,3.3){\begin{scriptsize}$0$\end{scriptsize}};
\node[circle, draw, minimum size=.2cm, fill=red] at (1,3){};
\node  at (1.95,3.3){\begin{tiny}$ac-b-c^2$\end{tiny}};
\node[circle, draw, minimum size=.2cm, fill=red] at (2,3){};
\node at (3.2, 3.3){\begin{scriptsize}$1$\end{scriptsize}};
\node[circle, draw, minimum size=.2cm, fill=red] at (3,3){};
\node at (2.0,1.7){\begin{scriptsize}$a-2c$\end{scriptsize}};
\node[circle, draw, minimum size=.2cm, fill=red] at (2,2){};
\node[circle, draw, minimum size=.2cm, fill=green] at (2,1){};
\end{tikzpicture}
\end{center}
\caption{Path set for $\gamma$}
\end{figure}

\section{Pseudo-involutions}
We note that by the theory of Riordan arrays, if we have $ac-b-c^2=0$, then the corresponding Riordan array will be a pseudo-involution \cite{PI}. In particular, this will be the case when
$$(a,b,c)=(r+1,r,r).$$
\begin{example} We consider the elliptic curve
$$y^2-3xy-y=x^3-2x^2-2x.$$
This leads to the Riordan array $(\gamma(x), x\gamma(x))$ that begins
$$\left(
\begin{array}{ccccccc}
 1 & 0 & 0 & 0 & 0 & 0 & 0 \\
 -1 & 1 & 0 & 0 & 0 & 0 & 0 \\
 1 &-2 & 1 & 0 & 0 & 0 & 0 \\
 0 & 3 & -3 & 1 & 0 & 0 & 0 \\
-2 & -2 & 6 & -4 & 1 & 0 & 0 \\
 5 & -3 & -7 & 10 & -6 & 1 & 0 \\
-7 & 14 & 0 & -16 & 15 & -6 & 1 \\
\end{array}
\right).$$ This is pseudo-involution in the Riordan group, which means that
$$(\gamma, -x\gamma)^2=I,$$ the identity matrix (represented by $(1,x)$ in the Riordan group).
\end{example}
Similarly, when $b=0$ and $a=c$, we will again obtain a pseudo-involution in the Riordan group.
\begin{example}
We consider the case $b=0$ and $a=c=-1$.
Thus we consider the elliptic curve
$$y^2+xy-y=x^3+x.$$
We find that $\gamma(x)$ expands to give the sequence \seqnum{A023431} that begins
$$1, 1, 1, 2, 4, 7, 13, 26, 52, 104, 212,\ldots.$$
The matrix $(\gamma, -x \gamma)$, which begins
$$\left(
\begin{array}{ccccccc}
 1 & 0 & 0 & 0 & 0 & 0 & 0 \\
 1 & -1 & 0 & 0 & 0 & 0 & 0 \\
 1 & -2 & 1 & 0 & 0 & 0 & 0 \\
 2 & -3 & 3 & 1 & 0 & 0 & 0 \\
 4 & -6 & 6 & -4 & 1 & 0 & 0 \\
 7 & -13 & 13 & -10 & 5 & 1 & 0 \\
13 & -26 & 30 & -24 & 15 & -6 & 1 \\
\end{array}
\right),$$ is a pseudo-involution in the Riordan group.

The Hankel transform of $\gamma(x)$ in these cases will be the periodic sequence \seqnum{A010892} that begins
$$0, -1, -1, 0, 1, 1, 0, -1, -1, 0, 1, 1, 0, -1, -1, 0,\ldots.$$
\end{example}

\section{Conclusions}
The content of this note can be summarized in the two following diagrams.
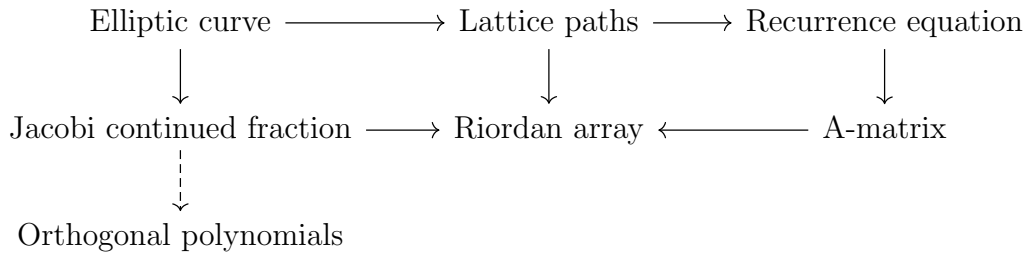
\begin{center}
\begin{tikzcd}
\text{Elliptic curve} \arrow[d] \arrow[r] & \text{Lattice paths} \arrow [d] \arrow[r] & \text{Recurrence equation} \arrow[d]\\
\text{Jacobi continued fraction} \arrow[d, dashed] \arrow[r] & \text{Riordan array}  & \text{A-matrix} \arrow[l]\\
\text{Orthogonal polynomials}
\end{tikzcd}
\end{center}
In this note, we have sought to show the close links between elliptic curves (containing the point $P=(0,0)$) and lattice paths, on the one hand, and elliptic division sequences (which in this case are Somos $4$ sequences) and the Hankel transforms of the lattice path enumerating sequences on the other. The group structure of the elliptic curve plays a strong role in both contexts, as the continued fraction expression of the generating function of the enumerator sequence is defined by the coordinates of the multiples of the point $P=(0,0)$ on the curve.
Note that we do not discuss the links to orthogonal polynomials here, but they have been explored elsewhere \cite{OP}. Solving the elliptic curve equation for $y$ and the reverting of a modified version of this solution play an important role. The modified version of the solution indicates that the dependence is not on the particular curve, but on its equivalence class with respect to the $j$-invariant. It is still unclear why the reversion process plays such an important role, so this suggests that further research into this could be fruitful.
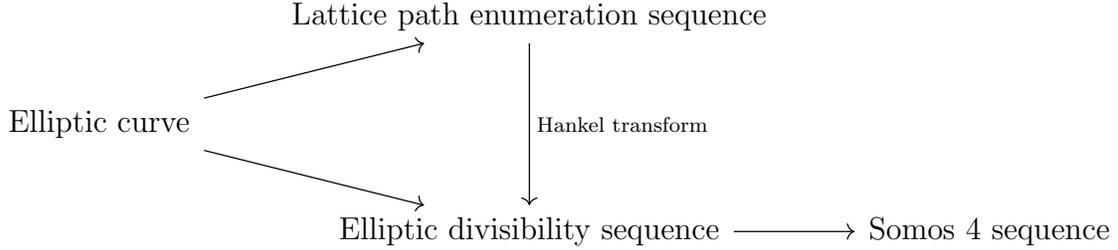
\begin{center}
\begin{tikzcd}
  & \text{Lattice path enumeration sequence} \arrow[dd,"\text{Hankel transform}"]\\
\text{Elliptic curve}  \arrow[rd] \arrow[ru] \\
 & \text{Elliptic divisibility sequence} \arrow[r] & \text{Somos $4$ sequence}
 \end{tikzcd}
 \end{center}

 \section{Appendix}
 The $A$-matrix approach \cite{He_M, Merlini} to the characterization of a Riordan array $(g(x), f(x))$ with general term $t_{n,k}$ leads to the following schema.
$$\left(
\begin{array}{ccccccc}
\vdots & \cdots & \cdots & \cdots & \cdots & \cdots & \cdots \\
x^2 & \cdots & t_{n-3,k-1} & t_{n-3,k} & t_{n-3,k+1} & t_{n-3,k+2} & \cdots \\
x & \cdots & t_{n-2,k-1} & t_{n-2,k} & t_{n-2,k+1} & t_{n-2,k+2} & \cdots \\
1 & \cdots & t_{n-1,k-1} & t_{n-1,k} & t_{n-1,k+1} & t_{n-1,k+2} & \cdots \\
1/x & \cdots & t_{n,k-1} & t_{n,k} & t_{n,k+1} & t_{n,k+2} & \cdots \\
1/x^2 & \cdots &  &  &  & t_{n+1,k+2} & \cdots \\
 &  &  &  &  &  &  \\
\cdots & \cdots & 1 & u & u^2 & u^3 & \cdots \\
\end{array}
\right).$$
This implies that if $t_{n,k}$ satisfies a recurrence such as
$$t_{n,k}=t_{n-1,k-1}+\gamma t_{n-2,k-1}+\alpha t_{n-1,k}+\beta t_{n-2,k}+\delta t_{n-2,k+1},$$ then $f(x)$ is defined by $u(x)=f(x)$ where we have
$$\frac{u}{x}=1+\gamma x+ \alpha u+ \beta ux+ \delta u^2x.$$
Solving for $u$, we find that
$$u=\frac{1-\alpha x- \beta x^2-\sqrt{1-2 \alpha x+(\alpha^2-2 \beta)x^2+2(\alpha \beta- 2 \delta)x^3+(\beta^2-4 \gamma \delta)x^4}}{2 \delta x^2}.$$
We let $g(x)=\frac{u(x)}{x}$. We then have
$$g(x)=\frac{1+\gamma x}{1-\alpha x- \beta x^2}C\left(\frac{\delta x^3(1+\gamma x)}{(1-\alpha x- \beta x^2)^2}\right).$$
The Hankel transform of the expansion of $g(x)$ is then a $(\delta^2, \delta^2(\alpha \gamma-\beta+\gamma^2))$ Somos $4$ sequence \cite{Chang, Yura}.

In this case, the $A$-matrix defining the corresponding Riordan array is given by
$$\left(\begin{array}{ccc}
\gamma & \beta & \delta\\
1 & \alpha & 0 \\
\end{array}\right).$$

\bigskip
\hrule

\noindent 2010 {\it Mathematics Subject Classification}:
Primary 05A15; Secondary 11G05, 14H52, 15B36, 11B37, 11B83.
\noindent \emph{Keywords:} Riordan array, lattice path, generating function, linear recurrence, elliptic curve, Somos sequence

\bigskip
\hrule
\bigskip
\noindent (Concerned with sequences
\seqnum{A000108},
\seqnum{A010892},
\seqnum{A023431}, and
\seqnum{A025243}
).


\begin{thebibliography}{99}


\bibitem{Book1} P. Barry, \emph{Riordan Arrays: a Primer}, Logic Press, 2017.

\bibitem{OP} P. Barry, Generalized Catalan recurrences, Riordan arrays, elliptic curves, and orthogonal polynomials, \url{https://arxiv.org/abs/1910.00875}.

\bibitem{Path} P. Barry, Notes on Riordan arrays and lattice paths, \url{https://arxiv.org/abs/2504.09719}.

\bibitem{PI} A. Burstein, L. W. Shapiro, Pseudo-involutions in the Riordan group, \url{https://arxiv.org/abs/2112.11595}.

\bibitem{Chang} Xiang-Ke Chang and Xing-Biao Hua, A conjecture based on Somos-$4$ sequence and its extension, \emph{Linear Algebra Appl.}, \textbf{436} (2012), 4285--4295.


\bibitem{He_M} Tian-Xiao He, Matrix characterizations of Riordan arrays, \emph{Linear Algebra Appl.}, \textbf{465} (2015), 15--42.

\bibitem{Hone} A. N. W. Hone, Elliptic Curves and Quadratic Recurrence Sequences, \emph{Bull. Lond. Math. Soc.}, \textbf{37} (2005), 161-171.

\bibitem{Layman} J. W. Layman, The Hankel transform and some of its properties, \emph{J. Integer Seq.}, \textbf{4} (2001),
\href{https://www.cs.uwaterloo.ca/journals/JIS/VOL4/LAYMAN/hankel.html} {Article 01.1.5}.

\bibitem{Merlini} D. Merlini, D. G. Rogers, R. Sprugnoli, M. C. Verri, On Some Alternative Characterizations of Riordan Arrays, \emph{Canad. J. Math.}, \textbf{49} (1997), 301--320.

\bibitem{Book2} L. Shapiro, R. Sprugnoli, P. Barry, G.-S. Cheon, T.-X. He, D. Merlini, and W. Wang, \emph{The Riordan Group and Applications}, Springer, 2022.

\bibitem{SGWW} L. W. Shapiro, S. Getu, W-J. Woan, and L.C. Woodson, The Riordan group, \emph{Discr. Appl. Math.}, \textbf{34} (1991),
 229--239.

\bibitem{Sil} J. H. Silverman, \emph{The Arithmetic of Elliptic Curves}, Springer.

\bibitem{SL1} N. J. A.~Sloane, \emph{The
On-Line Encyclopedia of Integer Sequences}. Published electronically
at \texttt{http://oeis.org}, 2025.

\bibitem{SL2} N. J. A.~Sloane, The On-Line Encyclopedia of Integer
Sequences, \emph{Notices Amer. Math. Soc.}, \textbf{50} (2003),  912--915.

\bibitem{Stan1} R. P. Stanley, \emph{Enumerative Combinatorics}, Volume 1, Cambridge University Press.

\bibitem{Stan2} R. P. Stanley, \emph{Enumberative Combinatorics}, Volume 2, Cambridge University Press.

\bibitem{StanC} R. P. Stanley, \emph{Catalan numbers}, Cambridge University Press.

\bibitem{Swart} C. Swart and A. N. W. Hone, Integrality and the Laurent phenomenon for Somos 4 sequences, \url{https://arxiv.org/abs/math/0508094}.

\bibitem{Wall} H. S. Wall, \emph{Analytic theory of continued fractions}, AMS Chelsea.

\bibitem{Wash} L. C. Washington, \emph{Elliptic Curves: Number Theory and Cryptography}, Routledge.

\bibitem{Yura} F. Yura, Hankel determinant solution for elliptic sequence, \emph{Linear Algebra Appl.}, \textbf{484} (2015), 27--45.

\end{thebibliography}
\end{document}